\documentclass{amsart}
\usepackage{lastpage}
\usepackage{yhmath,verbatim,mathabx}
\usepackage[all]{xy}

\DeclareFontFamily{U}{mathx}{\hyphenchar\font45}
\DeclareFontShape{U}{mathx}{m}{n}{
	<5> <6> <7> <8> <9> <10>
	<10.95> <12> <14.4> <17.28> <20.74> <24.88>
	mathx10
}{}
\DeclareSymbolFont{mathx}{U}{mathx}{m}{n}
\DeclareFontSubstitution{U}{mathx}{m}{n}
\DeclareMathAccent{\widecheck}{0}{mathx}{"71}
\DeclareMathAccent{\wideparen}{0}{mathx}{"75}

\newcommand{\bdis}{\begin{displaymath}}
\newcommand{\edis}{\end{displaymath}}
\newcommand{\be}{\begin{equation}}
\newcommand{\ee}{\end{equation}}
\newcommand{\mbb}{\mathbb}
\newcommand{\mcal}{\mathcal}

\newcommand{\vp}{\varphi}

\newcommand{\zf}{\zeta\left(\frac{1}{2}+it\right)}

\theoremstyle{definition}

\theoremstyle{remark}
\newtheorem{remark}[]{Remark}

\newtheorem*{mydef11}{{\bf Theorem 1}}

\newtheorem*{mydef12}{{\bf Theorem 2}}

\newtheorem*{mydef51}{{\bf Lemma 1}}

\newtheorem*{mydef52}{{\bf Lemma 2}}

\newtheorem*{mydef53}{{\bf Lemma 3}}

\newtheorem*{mydefECHF1}{{\bf Exact Complete Hybrid Formula 1}} 

\newtheorem*{mydefECHF2}{{\bf Exact Complete Hybrid Formula 2}} 

\newtheorem*{mydefETCHF}{{\bf Exact Ternary Complete Hybrid Formula }}

\numberwithin{equation}{section}

\begin{document}

\title{Jacob's ladders and new synergetic formula generating infinite set of $\zeta$-cooperative three-parametric invariants} 

\author{Jan Moser}

\address{Department of Mathematical Analysis and Numerical Mathematics, Comenius University, Mlynska Dolina M105, 842 48 Bratislava, SLOVAKIA}

\email{jan.mozer@fmph.uniba.sk}


\begin{abstract}
In this paper we obtain new $\zeta$-synergetic formula namely an exact secondary complete hybrid formula. This one is generated by some set of trigonometric and power functions together with the square of module of the Riemann's zeta-function on the critical line. By means of this new formula, we define a two-parametric continuum set of a three-parametric invariants. 
\end{abstract}
\maketitle

\section{Introduction}

\subsection{} 

This paper is based also on a set of new notions and methods we have introduced into the theory of the Riemann's zeta-function in our series of 47 papers concerning Jacob's ladders. These can be found in arXiv[math.CA] starting with the paper \cite{1}, (2009). 

Here we use especially the following notions: Jacob's ladder, $\zeta$-disconnected sets generated by the Jacob's ladder, algorithm for generating the $\zeta$-factorization formulas, complete hybrid formula (exact and asymptotic). Short survey of these notions can be found in our papers \cite{5}, \cite{8}. 

Next, we have introduced also the following notions (see \cite{8}): secondary crossbreeding, secondary complete hybrid formula and secondary $\zeta$-family of basic elements. 

\subsection{}  

In this paper we obtain some new results in the direction mentioned in 1.1. For example, we obtain as a simple variant of complete result, the following one:  the set 
\bdis 
\begin{split}
& \{ \sin^2t, \cos^2t, (t-\pi L)^{1/3}, (t-\pi L)^{1/5}\}, \\ 
& t\in [\pi L,\pi L+U],\ U\in (0,\pi/2),\ L\in\mbb{N} 
\end{split}
\edis 
of elementary functions generates the following exact secondary complete hybrid formula 
\be \label{1.1} 
\begin{split}
& \left( \frac{\alpha_0^{3,k_2}-\pi L}{\alpha_0^{4,k_2}-\pi L}\right)^{1/2}
\left\{\prod_{r=1}^{k_2}\frac{|\tilde{Z}(\alpha_r^{2,k_2})|^2|\tilde{Z}(\alpha_r^{3,k_2})|^3}{|\tilde{Z}(\alpha_r^{4,k_2})|^5}\right\}\cos^2(\alpha_0^{2,k_2})+\\ 
& + 
\left( \frac{\alpha_0^{3,k_1}-\pi L}{\alpha_0^{4,k_1}-\pi L}\right)^{1/2}
\left\{\prod_{r=1}^{k_1}\frac{|\tilde{Z}(\alpha_r^{1,k_1})|^2|\tilde{Z}(\alpha_r^{3,k_1})|^3}{|\tilde{Z}(\alpha_r^{4,k_1})|^5}\right\}\cos^2(\alpha_0^{1,k_1})= \\ 
& = \frac{81}{250}\sqrt{10},\ \forall\- L\geq L_0>0, 
\end{split}
\ee  
where 
\be \label{1.2}  
1\leq k_1,k_2,k_3,k_4\leq k_0, 
\ee  
and we fix arbitrary $k_0\in\mbb{N}$ and $L_0\in\mbb{N}$ being sufficiently big, 
\be \label{1.3} 
\tilde{Z}^2(t)=\frac{|\zf|^2}{\omega(t)},\ \omega(t)=\left\{1+\mcal{O}\left(\frac{\ln\ln t}{\ln t}\right)\right\}\ln t, 
\ee 
(see \cite{2}, (6.1), (6.7), (7.7), (7.8), (9.1)), next 
\be \label{1.4}  
\begin{split}
& \alpha_0^{1,k_1}, \alpha_0^{3,k_1}, \alpha_0^{4,k_1}, \alpha_0^{2,k_2}, \alpha_0^{3,k_2}, \alpha_0^{4,k_2}\in (\pi L,\pi L+U), \\ 
& \alpha_r^{1,k_1}, \alpha_r^{3,k_1}, \alpha_r^{4,k_1}\in 
(\overset{r}{\wideparen{\pi L}},\overset{r}{\wideparen{\pi L+U}}),\ r=1,\dots,k_1, \\ 
& \alpha_r^{2,k_2}, \alpha_r^{3,k_2}, \alpha_r^{4,k_2}\in 
(\overset{r}{\wideparen{\pi L}},\overset{r}{\wideparen{\pi L+U}}),\ r=1,\dots,k_2,
\end{split}
\ee 
and the segment 
\bdis 
[\overset{r}{\wideparen{\pi L}},\overset{r}{\wideparen{\pi L+U}}]
\edis  
is the $r$-th reverse iteration by means of the Jacob's ladder (see \cite{3}) of the basic segment 
\bdis 
[\pi L,\pi L+U]=[\overset{0}{\wideparen{\pi L}},\overset{0}{\wideparen{\pi L+U}}]. 
\edis  

\subsection{} 

Let us notice explicitly that the following holds true: consecutive components of the main $\zeta$-disconnected set 
\be \label{1.5} 
\Delta(\pi L,U,\bar{k})=\bigcup_{r=0}^{\bar{k}}
[\overset{r}{\wideparen{\pi L}},\overset{r}{\wideparen{\pi L+U}}]
\ee 
where 
\bdis 
\bar{k}=\max\{ k_1,k_2\},\ \bar{k}\leq k_0 
\edis 
are separated each from other by gigantic distances $\rho$: 
\be \label{1.6} 
\begin{split}
& \rho
\left\{
[\overset{r}{\wideparen{\pi L}},\overset{r}{\wideparen{\pi L+U}}]; 
[\overset{r+1}{\wideparen{\pi L}},\overset{r+1}{\wideparen{\pi L+U}}]
\right\}\sim (1-c)\pi(\pi L)\sim \\ 
& \sim \pi(1-c)\frac{L}{\ln L}\to\infty \ \mbox{as}\ L\to\infty,\ r=0,1,\dots \bar{k}-1, 
\end{split}
\ee 
where $c$ is the Euler's constant and $\pi(x)$ stands for the prime-counting function.  

\subsection{} 

Next (see (\ref{1.1}), (\ref{1.3})) the following asymptotic secondary hybrid formula 
\be \label{1.7} 
\begin{split}
& \left( \frac{\alpha_0^{3,k_2}-\pi L}{\alpha_0^{4,k_2}-\pi L}\right)^{1/2}
\left\{
\prod_{r=1}^{k_2}\frac{|\zeta(\frac 12+i\alpha_r^{2,k_2})|^2|\zeta(\frac 12+i\alpha_r^{3,k_3})|^3}
{|\zeta(\frac 12+i\alpha_r^{4,k_4})|^5}
\right\}\cos^2(\alpha_0^{2,k_2})+ \\ 
& + 
\left( \frac{\alpha_0^{3,k_1}-\pi L}{\alpha_0^{4,k_1}-\pi L}\right)^{1/2}
\left\{
\prod_{r=1}^{k_1}\frac{|\zeta(\frac 12+i\alpha_r^{1,k_1})|^2|\zeta(\frac 12+i\alpha_r^{3,k_1})|^3}
{|\zeta(\frac 12+i\alpha_r^{4,k_1})|^5}
\right\}\sin^2(\alpha_0^{1,k_1})\sim \\ 
& \sim \frac{81}{250}\sqrt{10}, \\ 
& 1\leq k_1,k_2\leq k_0,\ L\to\infty 
\end{split}
\ee  
holds true. 

\subsection{} 

Now, let us notice something about meaning of the the formula (\ref{1.1}) as well as the meaning of the formula (\ref{1.7}) introduced in the paper \cite{8}. Both mentioned formulas are synergetic ones follow as result of interactions controlled by the Jacob's ladder $\vp_1(t)$ between following continuum sets 
\bdis 
\left\{
\left|\zf\right|^2
\right\}, \{\sin^2t\}, \{\cos^2t\}, \{(t-\pi L)^{1/3}\}, \{(t-\pi L)^{1/5}\},\ t\geq L_0. 
\edis  
Let us call these interactions as the $\zeta$-chemical reaction as an analogue to the classical Belousov-Zhabotinski chemical oscillations. 

\begin{remark}
The result of the above mentioned $\zeta$-chemical reaction (the $\zeta$-chemical compound) is the synergetic formula (see (\ref{1.1}) and (\ref{1.7})) that contains a finite set of the arguments (see (\ref{1.4}) - for every fixed and admissible $U,L,k_1,k_2$, of course) that are distributed within disconnected set (\ref{1.5}) which components are separated by gigantic distances according to (\ref{1.6}). 
\end{remark} 

\subsection{} 

Let us notice explicitly about another result that is contained in the formula (\ref{1.1}). 

\begin{remark}
Namely, an invariant relatively to arbitrary sampling is defined by the left-hand side of the formula (\ref{1.1}): 
\begin{itemize}
	\item[(a)] of the continuous parameter $U\in(0,\pi/2)$ 
	\item[(b)] of the discrete parameter $L\geq L_0$ 
	\item[(c)] of the subset 
	\bdis 
	\{k_1,k_2\}\subset \{1,\dots,k_0\}
	\edis 
\end{itemize} 
\end{remark}

\begin{remark}
It is clear that the synergetic (cooperative) formula (\ref{1.1}) as well as (\ref{1.7}) represents completely new type of result in the theory of the Riemann's zeta-function and, simultaneously, in the theory of real continuous functions. 
\end{remark}

\section{The first exact complete hybrid formula} 

\subsection{} 

By making use of our algorithm for generating the $\zeta$-factorization formulas (see \cite{5}, (3.1) -- (3.11), comp. \cite{4}) we obtain the following result (comp. \cite{6}, (2.1) -- (2.9)): 

\begin{mydef51}
For the function 
\be \label{2.1} 
f_1(t)=\sin^2t\in \tilde{C}_0[\pi L,\pi L+U],\ U\in (0,\pi/2) 
\ee  
there are vector-valued functions 
\be \label{2.2} 
(\alpha_0^{1,k_1},\alpha_1^{1,k_1},\dots,\alpha_{k_1}^{1,k_1},\beta_1^{k_1},\dots,\beta_{k_1}^{k_1}),\ 1\leq k_1\leq k_0 
\ee  
such that the following exact $\zeta$-factorization formula 
\be \label{2.3} 
\prod_{r=1}^{k_1}\frac{\tilde{Z}^2(\alpha_r^{1,k_1})}{\tilde{Z}^2(\beta_r^{k_1})}=
\frac 12\left( 1-\frac{\sin 2U}{2U}\right)\frac{1}{\sin^2(\alpha_0^{1,k_1})},\ 
\forall\- L\geq L_0>0 
\ee  
holds true, (comp. \cite{5}, (3.7)), where 
\be \label{2.4} 
\begin{split}
& \alpha_r^{1,k_1}=\alpha_r(U,L,k_1;f_1),\ r=0,1,\dots,k_1, \\ 
& \beta_r^{k_1}=\beta_r(U,L,k_1),\ r=1,\dots,k_1, \\ 
& \alpha_0^{1,k_1}\in (\pi L,\pi L+U),\ 
\alpha_r^{1,k_1},\beta_r^{k_1}\in 
(\overset{r}{\wideparen{\pi L}},\overset{r}{\wideparen{\pi L+U}}),\ r=1,\dots,k_1. 
\end{split}
\ee 
\end{mydef51} 

\begin{mydef52}
For the function 
\be \label{2.5} 
f_2(t)=\cos^2t\in \tilde{C}_0[\pi L,\pi L+U],\ U\in (0,\pi/2) 
\ee  
there are vector-valued functions 
\be \label{2.6} 
(\alpha_0^{k_2},\alpha_1^{2,k_2},\dots,\alpha_{k_2}^{2,k_2},\beta_1^{k_2},\dots,\beta_{k_2}^{k_2}),\ 1\leq k_2\leq k_0 
\ee  
such that the following exact $\zeta$-factorization formula 
\be \label{2.7} 
\prod_{r=1}^{k_2}\frac{\tilde{Z}^2(\alpha_r^{2,k_2})}{\tilde{Z}^2(\beta_r^{k_2})}=
\frac 12\left( 1+\frac{\sin 2U}{2U}\right)\frac{1}{\cos^2(\alpha_0^{2,k_2})},\ 
\forall\- L\geq L_0>0 
\ee  
holds true, where 
\be \label{2.8} 
\begin{split}
	& \alpha_r^{2,k_2}=\alpha_r(U,L,k_2;f_2),\ r=0,1,\dots,k_2, \\ 
	& \beta_r^{k_2}=\beta_r(U,L,k_2),\ r=1,\dots,k_2, \\ 
	& \alpha_0^{2,k_2}\in (\pi L,\pi L+U),\ 
	\alpha_r^{2,k_2},\beta_r^{k_2}\in 
	(\overset{r}{\wideparen{\pi L}},\overset{r}{\wideparen{\pi L+U}}),\ r=1,\dots,k_2. 
\end{split}
\ee 
\end{mydef52} 

\subsection{} 

Now, we obtain by the crossbreeding (see \cite{6}, (2.9)) between $\zeta$-factorization formulae (\ref{2.3}) and (\ref{2.7}) the following 

\begin{mydefECHF1}
\be \label{2.9} 
\begin{split}
& 
\left\{
\prod_{r=1}^{k_2}\frac{\tilde{Z}^2(\alpha_r^{2,k_2})}{\tilde{Z}^2(\beta_r^{k_2})}
\right\}\cos^2(\alpha_0^{2,k_2})+
\left\{
\prod_{r=1}^{k_1}\frac{\tilde{Z}^2(\alpha_r^{1,k_1})}{\tilde{Z}^2(\beta_r^{k_1})}
\right\}\sin^2(\alpha_0^{1,k_1})=1, \\ 
& \forall\- L\geq L_0>0,\ 1\leq k_1,k_2\leq k_0. 
\end{split}
\ee 
\end{mydefECHF1} 

\section{The second exact complete hybrid formula} 

\subsection{} 

Further, we have the following exact $\zeta$-factorization formula (see \cite{6}, (3.1) -- (3.4), $L\to \pi L$). 

\begin{mydef53}
For the function 
\be \label{3.1} 
\bar{f}_\Delta(t)=\bar{f}(t;L,\Delta)=(t-\pi L)^\Delta\in \tilde{C}_0[\pi L,\pi L+U],\ U\in (0,\pi/2),\ \Delta>0 
\ee  
there are vector-valued functions 
\be \label{3.2} 
(\alpha_0^{\Delta,\bar{k}_\Delta},\alpha_1^{\Delta,\bar{k}_\Delta},\dots,\alpha_{\bar{k}_\Delta}^{\Delta,\bar{k}_\Delta},\beta_1^{\bar{k}_\Delta},\dots,\beta_{\bar{k}_\Delta}^{\bar{k}_\Delta}),\ 1\leq \bar{k}_\Delta\leq k_0,\ \bar{k}_\Delta\in\mbb{N}  
\ee  
such that the following exact $\zeta$-factorization formula 
\be \label{3.3} 
\prod_{r=1}^{\bar{k}_\Delta}\frac{\tilde{Z}^2(\alpha_r^{\Delta,\bar{k}_\Delta})}{\tilde{Z}^2(\beta_r^{\bar{k}_\Delta})}=
\frac{1}{1+\Delta}\left(\frac{U}{\alpha_0^{\Delta,\bar{k}_\Delta}-\pi L}\right)^\Delta,\ 
\forall\- L\geq L_0>0 
\ee  
holds true, where 
\be \label{3.4} 
\begin{split}
	& \alpha_r^{\Delta,\bar{k}_\Delta}=\alpha_r(U,L,\bar{k}_\Delta;\bar{f}_\Delta),\ r=0,1,\dots,\bar{k}_\Delta, \\ 
	& \beta_r^{\bar{k}_\Delta}=\beta_r(U,L,\bar{k}_\Delta),\ r=1,\dots,\bar{k}_\Delta, \\ 
	& \alpha_0^{\Delta,\bar{k}_\Delta}\in (\pi L,\pi L+U),\ 
	\alpha_r^{\Delta,\bar{k}_\Delta},\beta_r^{\bar{k}_\Delta}\in 
	(\overset{r}{\wideparen{\pi L}},\overset{r}{\wideparen{\pi L+U}}),\ r=1,\dots,\bar{k}_\Delta. 
\end{split}
\ee 
\end{mydef53} 

\subsection{} 

Next, we have the following set of two exact $\zeta$-factorization formulae (see (\ref{3.3})) 
\be \label{3.5} 
\begin{split}
& (1+\Delta_3)(\alpha_0^{3,k_3}-\pi L)^{\Delta_3}\prod_{r=1}^{k_3}\frac{\tilde{Z}^2(\alpha_r^{3,k_3})}{\tilde{Z}^2(\beta_r^{k_3})}=U^{\Delta_3}, \\ 
& (1+\Delta_4)(\alpha_0^{4,k_4}-\pi L)^{\Delta_4}\prod_{r=1}^{k_4}\frac{\tilde{Z}^2(\alpha_r^{4,k_4})}{\tilde{Z}^2(\beta_r^{k_4})}=U^{\Delta_4}, \\ 
& \Delta_3\not=\Delta_4, 
\end{split}
\ee  
where 
\be \label{3.6} 
\alpha_r^{\Delta_3,\bar{k}(\Delta_3)}=\alpha_r^{3,k_3},\ \dots 
\ee 
for the brevity. Now, we obtain, by crossbreeding on the set (\ref{3.5}) (i.e. by elimination of the external variable $U$ in (\ref{3.5})) the following 

\begin{mydefECHF2}
\be \label{3.7} 
\begin{split}
& (1+\Delta_3)^{1/\Delta_3}(\alpha_0^{3,k_3}-\pi L)
\left\{
\prod_{r=1}^{k_3}\frac{\tilde{Z}^2(\alpha_r^{3,k_3})}{\tilde{Z}^2(\beta_r^{k_3})}
\right\}^{1/\Delta_3}= \\ 
& = (1+\Delta_4)^{1/\Delta_4}(\alpha_0^{4,k_4}-\pi L)
\left\{
\prod_{r=1}^{k_4}\frac{\tilde{Z}^2(\alpha_r^{4,k_4})}{\tilde{Z}^2(\beta_r^{k_4})}
\right\}^{1/\Delta_4}, \\ 
& \forall\- L\geq L_0>0,\ \Delta_3,\Delta_4>0,\ \Delta_3\not=\Delta_4,\ 1\leq k_3,k_4\leq k_0. 
\end{split}
\ee 
\end{mydefECHF2} 

\section{The first variant of secondary exact complete hybrid formula} 

Now, we make use of our operation of the the secondary crossbreeding (see \cite{8}) on the set of two exact complete hybrid formulae (\ref{2.9}) and (\ref{3.7}) as follows. 

\subsection{} 

Firstly, we put into (\ref{3.7}) 
\be \label{4.1} 
k_3=k_4=k;\ 1\leq k\leq k_0 
\ee  
that gives us the following formula 
\be \label{4.2} 
\begin{split}
& \prod_{r=1}^k\tilde{Z}^2(\beta_r^k)=\\ 
& = 
\left[\frac{(1+\Delta_3)^{1/\Delta_3}}{(1+\Delta_4)^{1/\Delta_4}}\right]^{\frac{\Delta_3\Delta_4}{\Delta_4-\Delta_3}}
\left( \frac{\alpha_0^{3,k}-\pi L}{\alpha_0^{4,k}-\pi L}\right)^{\frac{\Delta_3\Delta_4}{\Delta_4-\Delta_3}}
\left\{
\prod_{r=1}^k\tilde{Z}^2(\alpha_r^{3,k})
\right\}^{\frac{\Delta_4}{\Delta_4-\Delta_3}}\times \\ 
& \times 
\left\{
\prod_{r=1}^k\tilde{Z}^2(\alpha_r^{4,k})
\right\}^{-\frac{\Delta_3}{\Delta_4-\Delta_3}}. 
\end{split}
\ee  

\subsection{} 

Secondly, we put consecutively in (\ref{4.2}) $k=k_1,k_2$ and the corresponding results we substitute into the formula (\ref{2.9}). Hence, we obtain the following 

\begin{mydef11} 
The set 
\be \label{4.3} 
\begin{split}
& \{\sin^2t,\cos^2t,(t-\pi L)^{\Delta_3},(t-\pi L)^{\Delta_4}\}, \\ 
& t\in [\pi L,\pi L+U],\ U\in (0,\pi/2),\ \Delta_3,\Delta_4>0, \Delta_3\not=\Delta_4 
\end{split}
\ee  
of real continuous functions generates the following secondary complete hybrid formula 
\be \label{4.4}
\begin{split}
& \frac
{\prod_{r=1}^{k_2}\tilde{Z}^2(\alpha_r^{2,k_2})[\tilde{Z}^2(\alpha_r^{3,k_2})]^{\frac{\Delta_4}{\Delta_3-\Delta_4}}[\tilde{Z}^2(\alpha_r^{4,k_2})]^{\frac{-\Delta_3}{\Delta_3-\Delta_4}}}
{\left( \frac{\alpha_0^{4,k_2}-\pi L}{\alpha_0^{3,k_2}-\pi L}\right)^{\frac{\Delta_3\Delta_4}{\Delta_3-\Delta_4}}}\cos^2(\alpha_0^{2,k_2})+ \\ 
& + 
\frac
{\prod_{r=1}^{k_1}\tilde{Z}^2(\alpha_r^{1,k_1})[\tilde{Z}^2(\alpha_r^{3,k_1})]^{\frac{\Delta_4}{\Delta_3-\Delta_4}}[\tilde{Z}^2(\alpha_r^{4,k_1})]^{\frac{-\Delta_3}{\Delta_3-\Delta_4}}}
{\left( \frac{\alpha_0^{4,k_1}-\pi L}{\alpha_0^{3,k_1}-\pi L}\right)^{\frac{\Delta_3\Delta_4}{\Delta_3-\Delta_4}}}\sin^2(\alpha_0^{1,k_1})= \\ 
& =
\left[
\frac
{(1+\Delta_4)^{1/\Delta_4}}
{(1+\Delta_3)^{1/\Delta_3}}
\right]^{\frac{\Delta_3\Delta_4}{\Delta_3-\Delta_4}},\ 
\forall\- L\geq L_0>0,\ 1\leq k_1,k_2\leq k_0 
\end{split}
\ee  
for arbitrary fixed $k_0\in\mbb{N}$ and sufficiently big $L_0\in\mbb{N}$. 
\end{mydef11}  

\begin{remark}
We list here main properties of the formula (\ref{4.4}): 
\begin{itemize}
	\item[(A)] it is the synergetic (cooperative) formula (see sect. 1.5 and paper \cite{8}), 
	\item[(B)] continuum set of the invariants relatively to arbitrary sampling of: 
	\begin{itemize}
		\item[(a)] of the continuum parameter $U\in(0,\pi/2)$, 
		\item[(b)] of discrete parameter $L\in [L_0,\infty)$, 
		\item[(c)] of the subset (as the set-parameter) 
		\bdis 
		\{k_1,k_2\}\subset \{1,\dots,k_0\}, 
		\edis 
	\end{itemize}
	is defined by its left-hand side for arbitrary 
	\bdis 
	\Delta_3,\Delta_4>0,\ \Delta_3\not=\Delta_4. 
	\edis 
	Consequently, we have defined two-parametric continuum set of three-parametric invariants. 
\end{itemize}
\end{remark} 

\begin{remark}
Formula (\ref{1.1}) corresponds to the values 
\bdis 
\Delta_3=\frac 13,\ \Delta_4=\frac 15
\edis 
in (\ref{4.4}). 
\end{remark} 

\section{Second variant of secondary exact complete hybrid formula}  

Here, we make again use of the secondary crossbreeding on the set of the formulae (\ref{2.9}) and (\ref{3.7}) however in the reverse direction as follows. 

\subsection{} 

First, we put in (\ref{2.9}) 
\be \label{5.1} 
k_1=k_2=k;\ 1\leq k\leq k_0 
\ee  
that gives us in result 
\be \label{5.2} 
\prod_{r=1}^k\tilde{Z}^2(\beta_r^k)=
\left\{
\prod_{r=1}^k\tilde{Z}^2(\alpha_r^{2,k})
\right\}\cos^2(\alpha_0^{2,k})+
\left\{
\prod_{r=1}^k\tilde{Z}^2(\alpha_r^{1,k})
\right\}\sin^2(\alpha_0^{1,k}). 
\ee  

\subsection{} 

Secondly, we put consecutively in (\ref{5.2}) 
\bdis 
k=k_3,k_4 
\edis  
and we put the corresponding result into (\ref{3.7}). The final result is listed bellow as theorem. 

\begin{mydef12}
The set 
\be \label{5.3} 
\begin{split}
	& \{\sin^2t,\cos^2t,(t-\pi L)^{\Delta_3},(t-\pi L)^{\Delta_4}\}, \\ 
	& t\in [\pi L,\pi L+U],\ U\in (0,\pi/2),\ \Delta_3,\Delta_4>0, \Delta_3\not=\Delta_4 
\end{split}
\ee  
generates the following secondary complete hybrid formula 
\be \label{5.4} 
\begin{split}
& \frac{\alpha_0^{3,k_3}-\pi L}{\alpha_0^{3,k_3}-\pi L}\times \\ 
& \times 
\left( 
\frac
{\prod_{r=1}^{k_3}\tilde{Z}^2(\alpha_r^{3,k_3})}
{[\prod_{r=1}^{k_3}\tilde{Z}^2(\alpha_r^{2,k_3})]\cos^2(\alpha_0^{2,k_3})+
[\prod_{r=1}^{k_3}\tilde{Z}^2(\alpha_r^{1,k_3})]\sin^2(\alpha_0^{1,k_3})}
\right)^{\frac{1}{\Delta_3}}\times \\ 
& \times 
\left( 
\frac
{\prod_{r=1}^{k_4}\tilde{Z}^2(\alpha_r^{4,k_4})}
{[\prod_{r=1}^{k_4}\tilde{Z}^2(\alpha_r^{2,k_4})]\cos^2(\alpha_0^{2,k_4})+
	[\prod_{r=1}^{k_4}\tilde{Z}^2(\alpha_r^{1,k_4})]\sin^2(\alpha_0^{1,k_4})}
\right)^{-\frac{1}{\Delta_4}}= \\ 
& = 
\frac
{(1+\Delta_4)^{1/\Delta_4}}
{(1+\Delta_3)^{1/\Delta_3}},\ \forall\- L\geq L_0>0,\ 1\leq k_3,k_4\leq k_0 . 
\end{split}
\ee 
\end{mydef12} 

\begin{remark}
The case 
\bdis 
\Delta_3=\Delta_4 
\edis  
is excluded (see (\ref{5.3})) since this one implies trivially that $1=1$. 
\end{remark} 

\begin{remark}
Main properties of the formula (\ref{5.4}) are characterized by the Remark 4. 
\end{remark} 

\section{On notion of exact ternary complete hybrid formula} 

Now, we are able to proceed in our hierarchy of types of crossbreeding. Namely, we can introduce notion of the exact  ternary complete hybrid formula. Of course, it is sufficient to give an example. 

We have obtained the set of two exact secondary hybrid formulae (\ref{4.4}), (\ref{5.4}). Every of these formulae contains the functions (this is the point to be generalized) 
\bdis 
\frac{(1+\Delta_4)^{1/\Delta_4}}{(1+\Delta_3)^{1/\Delta_3}},\ \Delta_3,\Delta_4>0,\ \Delta_3\not=\Delta_4. 
\edis  
The elimination of this function from the set $\{(\ref{4.4}), (\ref{5.4})\}$ (that is the ternary crossbreeding) gives the following 

\begin{mydefETCHF}
\be \label{6.1} 
\begin{split}
& \frac
{\prod_{r=1}^{k_2}\tilde{Z}^2(\alpha_r^{2,k_2})[\tilde{Z}^2(\alpha_r^{3,k_2})]^{\frac{\Delta_4}{\Delta_3-\Delta_4}}[\tilde{Z}^2(\alpha_r^{4,k_2})]^{\frac{-\Delta_3}{\Delta_3-\Delta_4}}}
{\left( \frac{\alpha_0^{4,k_2}-\pi L}{\alpha_0^{3,k_2}-\pi L}\right)^{\frac{\Delta_3\Delta_4}{\Delta_3-\Delta_4}}}\cos^2(\alpha_0^{2,k_2})+ \\ 
& + 
\frac
{\prod_{r=1}^{k_1}\tilde{Z}^2(\alpha_r^{1,k_1})[\tilde{Z}^2(\alpha_r^{3,k_1})]^{\frac{\Delta_4}{\Delta_3-\Delta_4}}[\tilde{Z}^2(\alpha_r^{4,k_1})]^{\frac{-\Delta_3}{\Delta_3-\Delta_4}}}
{\left( \frac{\alpha_0^{4,k_1}-\pi L}{\alpha_0^{3,k_1}-\pi L}\right)^{\frac{\Delta_3\Delta_4}{\Delta_3-\Delta_4}}}\sin^2(\alpha_0^{1,k_1})= \\ 
& \left(\frac{\alpha_0^{3,k_3}-\pi L}{\alpha_0^{3,k_3}-\pi L}\right)^{\frac{\Delta_3\Delta_4}{\Delta_3-\Delta_4}}\times \\ 
& \times 
\left( 
\frac
{\prod_{r=1}^{k_3}\tilde{Z}^2(\alpha_r^{3,k_3})}
{[\prod_{r=1}^{k_3}\tilde{Z}^2(\alpha_r^{2,k_3})]\cos^2(\alpha_0^{2,k_3})+
	[\prod_{r=1}^{k_3}\tilde{Z}^2(\alpha_r^{1,k_3})]\sin^2(\alpha_0^{1,k_3})}
\right)^{\frac{\Delta_4}{\Delta_3-\Delta_4}}\times \\ 
& \times 
\left( 
\frac
{\prod_{r=1}^{k_4}\tilde{Z}^2(\alpha_r^{4,k_4})}
{[\prod_{r=1}^{k_4}\tilde{Z}^2(\alpha_r^{2,k_4})]\cos^2(\alpha_0^{2,k_4})+
	[\prod_{r=1}^{k_4}\tilde{Z}^2(\alpha_r^{1,k_4})]\sin^2(\alpha_0^{1,k_4})}
\right)^{-\frac{\Delta_3}{\Delta_3-\Delta_4}}, \\ 
& \forall\- L\geq L_0>0,\ 1\leq k_1,k_2,k_3,k_4\leq k_0 .
\end{split}
\ee 
\end{mydefETCHF}

\end{document}